# The symplectic Thom conjecture

By Peter Ozsváth and Zoltán Szabó*

**Abstract**

In this paper, we demonstrate a relation among Seiberg-Witten invariants which arises from embedded surfaces in four-manifolds whose self-intersection number is negative. These relations, together with Taubes' basic theorems on the Seiberg-Witten invariants of symplectic manifolds, are then used to prove the symplectic Thom conjecture: a symplectic surface in a symplectic four-manifold is genus-minimizing in its homology class. Another corollary of the relations is a general adjunction inequality for embedded surfaces of negative self-intersection in four-manifolds.

## 1. Introduction

An old conjecture attributed to Thom states that the smooth holomorphic curves of degree $d$ in $\mathbb{C}P^2$ are genus-minimizing in their homology class. In their seminal paper, Kronheimer and Mrowka [12] showed how techniques from gauge theory can be brought to bear on questions of this kind: they showed that a "generalized Thom conjecture" holds for algebraic curves with nonnegative self-intersection number in a wide class of Kähler surfaces (which excludes $\mathbb{C}P^2$). With the advance of the Seiberg-Witten equations [29], Kronheimer-Mrowka [11] and Morgan-Szabó-Taubes [18] proved the Thom conjecture for curves with nonnegative self-intersection in any Kähler surface. Indeed, in light of Taubes' ground-breaking results [26], [27], the results readily generalized to the symplectic context; see [18], [14]. However, these proofs all hinged on the assumption that the surface has nonnegative self-intersection; the case where the self-intersection is negative remained elusive, save for a result by Fintushel and Stern for immersed spheres [6]; see also [25]. Our goal here is to prove the symplectic Thom conjecture in its complete generality (cf. Kirby's problem list [10, p. 326]):

*The first author was partially supported by NSF grant number DMS 9304580. The second author was partially supported by NSF grant number DMS 970435 and a Sloan Fellowship.

THEOREM 1.1 (symplectic Thom conjecture). *An embedded symplectic surface in a closed, symplectic four-manifold is genus-minimizing in its homology class.*

The following special case is of interest in its own right.

COROLLARY 1.2 (Kähler case). *An embedded holomorphic curve in a Kähler surface is genus-minimizing in its homology class.*

The proof is based on a relation among Seiberg-Witten invariants. Before stating the relation, we set up notation. Let $X$ be a closed, smooth four-manifold equipped with an orientation for which $b_2^+(X) > 0$ (here, $b_2^+(X)$ is the dimension of a maximal positive-definite linear subspace $H^+(X;\mathbb{R})$ of the intersection pairing on $H^2(X;\mathbb{R})$) and an orientation for $H^1(X;\mathbb{R}) \oplus H^+(X;\mathbb{R})$ (the latter is called a *homology orientation*). Given such a four-manifold, together with a $\mathrm{Spin}_{\mathbb{C}}$ structure $\mathfrak{s}$, the Seiberg-Witten invariants form an integer-valued function

$$\mathrm{SW}_{X,\mathfrak{s}} \colon \mathbb{A}(X) \longrightarrow \mathbb{Z},$$

where and $\mathbb{A}(X)$ denotes the graded algebra obtained by tensoring the exterior algebra on $H_1(X)$ (graded so that $H_1(X)$ has grading one) with the polynomial algebra $\mathbb{Z}[U]$ on a single two-dimensional generator. (We drop the four-manifold $X$ from the notation when it is clear from the context.) Each $\mathrm{Spin}_{\mathbb{C}}$ structure is specified by a pair of unitary $\mathbb{C}^2$ bundles $W^+$ and $W^-$ (the bundles of spinors), together with a Clifford action

$$\rho \colon T^*X \otimes W^+ \longrightarrow W^-.$$

Given a $\mathrm{Spin}_{\mathbb{C}}$ structure $\mathfrak{s}$ over $X$, the Seiberg-Witten invariant of $\mathfrak{s}$ is defined by "counting" the number of spin-connections $A$ on $W^+$ and spinors $\Phi \in \Gamma(X, W^+)$, up to gauge, which satisfy the Seiberg-Witten equations

$$\begin{align}
\rho_{\Lambda^+}(\mathrm{Tr} F_A^+) &= i(\Phi \otimes \Phi^*)_\circ + \rho_{\Lambda^+}(i\eta) \tag{1}\\
\slashed{D}_A^+ \Phi &= 0, \tag{2}
\end{align}$$

where $\eta$ is some fixed, real self-dual two-form. Here, $\mathrm{Tr} F_A^+$ denotes the trace of the self-dual part of the curvature form of $A$ (or, equivalently, the self-dual part of the Chern-Weil representative for $c_1(W^+)$ induced from $A$), $\rho_{\Lambda^+}$ denotes the endomorphism of $W^+$ induced from the Clifford action of self-dual two-forms, which in turn is induced from the Clifford action of one-forms, $(\Phi \otimes \Phi^*)_\circ$ denotes the endomorphism of $W^+$ which maps any spinor $\psi \in \Gamma(X, W^+)$ to

$$(\Phi \otimes \Phi^*)_\circ \psi = \langle \psi, \Phi \rangle \Phi - \frac{|\Phi|^2}{2} \psi,$$

and $\slashed{D}_A^+$ denotes the $\mathrm{Spin}_{\mathbb{C}}$ Dirac operator taking $W^+$ to $W^-$, coupled to the connection $A$. In the definition of the invariant, solutions are to be counted in the following sense. If $\mathcal{M}(X,\mathfrak{s})$ denotes the moduli space of solutions to equations (1) and (2) modulo gauge, then

$$\mathrm{SW}_{X,\mathfrak{s}}(a) = \langle \mu(a), [\mathcal{M}(X,\mathfrak{s})]\rangle,$$

where $[\mathcal{M}(X,\mathfrak{s})]$ denotes the fundamental class for the moduli space induced from the homology orientation, and

$$\mu\colon \mathbb{A}(X) \longrightarrow H^*(\mathcal{M}(X,\mathfrak{s}); \mathbb{Z})$$

denotes the map given by

$$\mu(x) = c_1(\mathcal{L})/x,$$

for $x \in H_*(X;\mathbb{Z})$ (for this purpose, we consider $U$ to be a generator of $H_0(X;\mathbb{Z})$), where $\mathcal{L}$ is the universal bundle over $X \times \mathcal{M}(X,\mathfrak{s})$. The pairing defining the invariant is nonzero only on those homogeneous elements $a$ whose degree is $d(\mathfrak{s})$, the expected dimension of the moduli space. The Atiyah-Singer index theorem allows one to express this dimension in terms of the homotopy type of $X$ and the first Chern class of $\mathfrak{s}$, $c_1(\mathfrak{s})$ (which is defined to be $c_1(W^+)$):

$$(3) \qquad d(\mathfrak{s}) = \frac{c_1(\mathfrak{s})^2 - (2\chi(X) + 3\sigma(X))}{4}$$

(here, $\chi(X)$ is the Euler characteristic of $X$ and $\sigma(X)$ is the signature of its intersection form on $H^2(X;\mathbb{R})$), or equivalently, $d(\mathfrak{s})$ is the Euler number of the bundle $W^+$.

The Seiberg-Witten invariant is a smooth invariant of the four-manifold $X$ (i.e. independent of the metric and perturbation $\eta$) when $b_2^+(X) > 1$. When $b_2^+(X) = 1$, then the invariant depends on the chamber, as follows. Let

$$\Omega(X) = \{x \in H^2(X;\mathbb{R}) \big| x^2 = 1\};$$

then $\Omega(X)$ has two components, and an orientation of $H^+(X;\mathbb{R})$ determines the positive component, denoted $\Omega^+(X)$. For a given $\mathrm{Spin}_{\mathbb{C}}$ structure $\mathfrak{s}$, the *wall determined by* $\mathfrak{s}$, denoted $\mathcal{W}_{\mathfrak{s}}$, is the set of $(\omega, t) \in \Omega^+(X) \times \mathbb{R}$ so that $2\pi\omega \cdot c_1(\mathfrak{s}) + t = 0$. The *chamber determined by* $\mathfrak{s}$ is a connected component of $(\Omega^+(X) \times \mathbb{R}) - \mathcal{W}_{\mathfrak{s}}$. There is a map, the *period map*, from the space of metrics and perturbations to the space $\Omega^+(X) \times \mathbb{R}$ defined by taking $g$ and $\eta$ to $\omega_g$ and $t = \int_X \omega_g \wedge \eta$, where $\omega_g$ is the unique harmonic, self-dual two-form $\omega_g \in \Omega^+(X)$. The Seiberg-Witten invariant of $\mathfrak{s}$ for $g$ and $\eta$ is well-defined if the corresponding period point does not lie on a wall, and it depends on $g$ and $\eta$ only through the chamber of the associated period point (see [11] and [14]). Let $\mathfrak{S}$ be a collection of $\mathrm{Spin}_{\mathbb{C}}$ structures, then a *common chamber for* $\mathfrak{S}$ is a connected component of

$$\Omega^+(X) \times \mathbb{R} - \bigcup_{\mathfrak{s} \in \mathfrak{S}} \mathcal{W}_{\mathfrak{s}}.$$

Given a cohomology class $c \in H^2(X;\mathbb{Z})$ of negative square, a common chamber for $\mathfrak{S}$ is called *perpendicular to $c$* if it contains a pair $(\omega, t)$, with $\omega$ perpendicular to $c$.

Given a $\mathrm{Spin}_{\mathbb{C}}$ structure $\mathfrak{s}$ determined by $(\rho, W^+, W^-)$, and a Hermitian line bundle $L$, there is a new $\mathrm{Spin}_{\mathbb{C}}$ structure corresponding to

$$(\rho, W^+ \otimes L, W^- \otimes L).$$

The isomorphism class of this new $\mathrm{Spin}_{\mathbb{C}}$ structure depends only on the first Chern class $c$ of $L$, so we will denote the new structure by $\mathfrak{s} + c$.

Let $\Sigma \subset X$ be a smoothly embedded, oriented, closed surface. We define the class $\xi(\Sigma) \in \mathbb{A}(X)$ by

$$\xi(\Sigma) = \prod_{i=1}^{g}(U - A_i \cdot B_i),$$

where $\{A_i, B_i\}_{i=1}^{g}$ are the images in $H_1(X;\mathbb{Z})$ of a standard symplectic basis for $H_1(\Sigma;\mathbb{Z})$, and the product $A_i \cdot B_i$ is taken in the algebra $\mathbb{A}(X)$. Of course, $\xi(\Sigma)$ depends on the orientation of $\Sigma$.

We can now state the relation.

THEOREM 1.3. *Let $X$ be a closed, smooth four-manifold with $b_2^+(X) > 0$ and $\Sigma \subset X$ a smoothly embedded, oriented, closed surface of genus $g > 0$ and negative self-intersection number*

$$[\Sigma] \cdot [\Sigma] = -n.$$

*If $b_2^+(X) > 1$, then for each $\mathrm{Spin}_{\mathbb{C}}$ structure $\mathfrak{s}$ with $d(\mathfrak{s}) \geq 0$ and*

$$|\langle c_1(\mathfrak{s}), [\Sigma]\rangle| \geq 2g + n,$$

*we have for each $a \in \mathbb{A}(X)$,*

(4) $$\mathrm{SW}_{\mathfrak{s}+\epsilon\mathrm{PD}(\Sigma)}(\xi(\epsilon\Sigma)U^m \cdot a) = \mathrm{SW}_{\mathfrak{s}}(a),$$

*where $\epsilon = \pm 1$ is the sign of $\langle c_1(\mathfrak{s}), [\Sigma]\rangle$, $2m = |\langle c_1(\mathfrak{s}), [\Sigma]\rangle| - 2g - n$, and $\mathrm{PD}(\Sigma)$ denotes the class Poincaré dual to $[\Sigma]$. Furthermore, if $b_2^+(X) = 1$, then the above relation holds in any common chamber for $\mathfrak{s}$ and $\mathfrak{s} + \epsilon\mathrm{PD}(\Sigma)$ which is perpendicular to $\mathrm{PD}(\Sigma)$.*

*Remark* 1.4. Let

$$\Sigma^{\perp} = \{(\omega, t) \in \Omega^+(X) \times \mathbb{R} \,\big|\, \omega \cdot \mathrm{PD}(\Sigma) = 0\}.$$

Note that

$$\Sigma^{\perp} \cap \mathcal{W}_{\mathfrak{s}} = \Sigma^{\perp} \cap \mathcal{W}_{\mathfrak{s}+\epsilon\mathrm{PD}(\Sigma)},$$

so that there are exactly two common chambers for $\mathfrak{s}$ and $\mathfrak{s} + \epsilon\mathrm{PD}(\Sigma)$ which are perpendicular to $\mathrm{PD}(\Sigma)$.

*Remark* 1.5. The above theorem holds also when $g = 0$. This was proved by Fintushel and Stern [6].

*Remark* 1.6. Some authors consider the Seiberg-Witten invariant as a function depending on characteristic cohomology elements $K \in H^2(X; \mathbb{Z})$,

$$\overline{\mathrm{SW}}_{X,K} \colon \mathbb{A}(X) \longrightarrow \mathbb{Z}$$

which is related to the more refined invariant used here by

$$\overline{\mathrm{SW}}_{X,K} = \sum_{\{\mathfrak{s} | c_1(\mathfrak{s}) = K\}} \mathrm{SW}_{X,\mathfrak{s}}.$$

Since $c_1(\mathfrak{s} + c) = c_1(\mathfrak{s}) + 2c$, the relation (4) implies a relation

$$\overline{\mathrm{SW}}_{X, K+2\epsilon \mathrm{PD}(\Sigma)}(\xi(\epsilon \Sigma) U^m \cdot a) = \overline{\mathrm{SW}}_{X,K}(a).$$

An important, and immediate, consequence of this relation is the following adjunction inequality for "basic classes" of four-manifolds of simple type $X$. A *basic class* is a $\mathrm{Spin}_{\mathbb{C}}$ structure $\mathfrak{s}$ whose Seiberg-Witten invariant does not vanish identically, and a four-manifold is said to be *of simple type* if all of the moduli spaces associated with its basic classes are zero-dimensional. More generally, a four-manifold is said to be *of type m* if the Seiberg-Witten invariants vanish for all $\mathrm{Spin}_{\mathbb{C}}$ structures $\mathfrak{s}$ with

$$d(\mathfrak{s}) \geq 2m.$$

COROLLARY 1.7 (adjunction inequality for negative self-intersections). *Let $X$ be a four-manifold of Seiberg-Witten simple type with $b_2^+(X) > 1$, and $\Sigma \subset X$ be a smoothly embedded, oriented, closed surface of genus $g(\Sigma) > 0$ and $\Sigma \cdot \Sigma < 0$. Then for all Seiberg-Witten basic classes $\mathfrak{s}$,*

(5) $$|\langle [\Sigma], c_1(\mathfrak{s}) \rangle| + [\Sigma] \cdot [\Sigma] \leq 2g(\Sigma) - 2.$$

*Proof.* If we had a basic class $\mathfrak{s}$ which violated the adjunction inequality, then the relation would guarantee that $\mathfrak{s} + \epsilon \mathrm{PD}(\Sigma)$ is also a basic class. But

$$d(\mathfrak{s} + \epsilon \mathrm{PD}(\Sigma)) = d(\mathfrak{s}) + \epsilon \langle c_1(\mathfrak{s}), [\Sigma] \rangle - n \geq d(\mathfrak{s}) + 2g > 0,$$

which violates the simple type assumption. □

*Remark* 1.8. A similar result for immersed spheres was proved by Fintushel and Stern in [6].

*Remark* 1.9. Note that the above argument in fact shows that the adjunction inequality (5) holds for all four-manifolds of type $g$.

*Remark* 1.10. We have the following analogous statement when $b_2^+(X) = 1$. Let $\mathfrak{S}_+$ denote the set of all $\mathrm{Spin}_\mathbb{C}$ structures $\mathfrak{s}$ with $d(\mathfrak{s}) \geq 0$. A common chamber for $\mathfrak{S}_+$ is called *of type m* if all its Seiberg-Witten invariants vanish for the $\mathrm{Spin}_\mathbb{C}$ structures $\mathfrak{s}$ with $d(\mathfrak{s}) \geq 2m$. The above proof also shows that the adjunction inequality holds in those common chambers for $\mathfrak{S}_+$ which are perpendicular to $\mathrm{PD}(\Sigma)$ and of type $g$.

*Remark* 1.11. By further extending the methods begun in this paper, one can obtain other adjunction inequalities for nonsimple type four-manifolds; see [21].

This paper is organized as follows. In Section 2, we prove the relation (Theorem 1.3), assuming results about the moduli spaces over a tubular neighborhood $N$ of $\Sigma$ (using a suitable connection on $TN$). In Section 3, we show how Theorem 1.1 follows from the relation. In Sections 4–6, we address the technical points assumed in Section 2. In Section 4 we show that the Seiberg-Witten invariant can be calculated using any connection on the tangent bundle, by modifying the usual compactness arguments (see Theorem 4.6). In Section 5, we construct a particularly convenient connection on $N$ for which the Dirac operator admits a holomorphic description. That section then concludes with some immediate consequences for moduli spaces over $N$. In Section 6, the holomorphic description is used to describe the obstruction bundles over the moduli spaces for $N$, completing the argument from Section 2. In Section 7, we give some examples of our results.

## 2. The relation

The goal of this section is to prove Theorem 1.3, assuming some technical facts which are proved in Sections 4–6. First, we show how to reduce the theorem to the following special case:

PROPOSITION 2.1. *Let $X$, $\Sigma$, $g$, and $n$ be as in Theorem 1.3. Assume moreover that $n \geq 2g$ and that the $\mathrm{Spin}_\mathbb{C}$ structure $\mathfrak{s}$ satisfies the condition*

$$\langle c_1(\mathfrak{s}), [\Sigma] \rangle = -2g - n.$$

*Then Relation* (4) *holds.*

The reduction involves the following blow-up formula of Fintushel and Stern:

THEOREM 2.2 (blow-up formula [6] and [22]). *Let $X$ be a smooth, closed four-manifold, and let $\widehat{X} = X \# \overline{\mathbb{C}P}^2$ denote its blow-up, with exceptional class $E \in H^2(\widehat{X}; \mathbb{Z})$. If $b_2^+(X) > 1$, then for each $\mathrm{Spin}_\mathbb{C}$ structure $\widehat{\mathfrak{s}}$ on $\widehat{X}$ with*

$d(\widehat{\mathfrak{s}}) \geq 0$, and each $a \in \mathbb{A}(X) \cong \mathbb{A}(\widehat{X})$, we have

$$\mathrm{SW}_{\widehat{X},\widehat{\mathfrak{s}}}(a) = \mathrm{SW}_{X,\mathfrak{s}}(U^m a),$$

where $\mathfrak{s}$ is the $\mathrm{Spin}_{\mathbb{C}}$ structure induced on $X$ obtained by restricting $\widehat{\mathfrak{s}}$, and $2m = d(\mathfrak{s}) - d(\widehat{\mathfrak{s}})$. If $b_2^+(X) = 1$, then for each chamber $\widehat{C}$ perpendicular to $E$, we have the above relation, where $\mathrm{SW}_{X,\mathfrak{s}}$ is calculated in the chamber on $X$ induced from $\widehat{C}$.

Now, we turn to the proof of Theorem 1.3 assuming Proposition 2.1, whose proof we will give afterwards.

*Proof of Theorem* 1.3. By reversing the orientation of $\Sigma$, it suffices to prove the relation when

$$\langle c_1(\mathfrak{s}), [\Sigma] \rangle \leq -n - 2g.$$

Let

$$m = \frac{1}{2}(-\langle c_1(\mathfrak{s}), [\Sigma] \rangle - n - 2g).$$

We blow up our manifold $\ell + m$ times to obtain a new manifold

$$\widehat{X} = X \# (\ell + m) \overline{\mathbb{C}P^2},$$

where $\ell$ is chosen so that

$$n + \ell + m \geq 2g.$$

Let $\widehat{\Sigma}$ be the proper transform of $\Sigma$, i.e. an embedded submanifold of the same genus for which

$$\mathrm{PD}(\widehat{\Sigma}) = \mathrm{PD}(\Sigma) - E_1 - \ldots - E_{\ell+m}.$$

Consider the $\mathrm{Spin}_{\mathbb{C}}$ structure $\widehat{\mathfrak{s}}$ on $\widehat{X}$ which extends $\mathfrak{s}$ on $X$ and whose first Chern class satisfies

$$c_1(\widehat{\mathfrak{s}}) = c_1(\mathfrak{s}) - E_1 - \ldots - E_\ell + E_{\ell+1} + \ldots + E_{\ell+m}.$$

Note that

$$\begin{aligned} [\widehat{\Sigma}] \cdot [\widehat{\Sigma}] &= -n - \ell - m, \\ \langle c_1(\widehat{\mathfrak{s}}), [\widehat{\Sigma}] \rangle &= -2g - n - \ell - m \\ d(\widehat{\mathfrak{s}}) &= d(\mathfrak{s}). \end{aligned}$$

We can now apply Proposition 2.1 to $\widehat{X}$, $\widehat{\Sigma}$, and $\widehat{\mathfrak{s}}$, to conclude that

(6) $$\mathrm{SW}_{\widehat{X},\widehat{\mathfrak{s}}}(a) = \mathrm{SW}_{\widehat{X},\widehat{\mathfrak{s}}-\mathrm{PD}[\widehat{\Sigma}]}(\xi(-\widehat{\Sigma}) \cdot a).$$

Since

$$d(\widehat{\mathfrak{s}} - \mathrm{PD}[\widehat{\Sigma}]) \geq d(\widehat{\mathfrak{s}}) = d(\mathfrak{s}) \geq 0,$$

and

$$(\widehat{\mathfrak{s}} - \mathrm{PD}[\widehat{\Sigma}])|_X = \mathfrak{s} - \mathrm{PD}[\Sigma],$$

the blow-up formula says that

(7) $$\mathrm{SW}_{\widehat{X},\widehat{\mathfrak{s}}-\mathrm{PD}[\widehat{\Sigma}]}(\xi(-\widehat{\Sigma}) \cdot a) = \mathrm{SW}_{X,\mathfrak{s}-\mathrm{PD}[\Sigma]}(\xi(-\Sigma)U^m \cdot a)$$

and

(8) $$\mathrm{SW}_{\widehat{X},\widehat{\mathfrak{s}}}(a) = \mathrm{SW}_{X,\mathfrak{s}}(a).$$

Putting together equations (6), (7), and (8), we prove the relation. In the case when $b_2^+(X) = 1$, a common chamber for $\mathfrak{s}$ and $\mathfrak{s} - \mathrm{PD}(\Sigma)$ which is orthogonal to $\Sigma$ gives rise to a common chamber for $\widehat{\mathfrak{s}}$ and $\widehat{\mathfrak{s}} - \mathrm{PD}(\widehat{\Sigma})$ which is orthogonal to $\widehat{\Sigma}$. Using this chamber, the above arguments go through. □

The rest of this section is devoted to the proof of Proposition 2.1. We assume $n \geq 2g$ and $\mathfrak{s}_0$ is a $\mathrm{Spin}_{\mathbb{C}}$ structure for which

(9) $$\langle c_1(\mathfrak{s}_0), [\Sigma] \rangle = -2g - n.$$

The proof involves stretching the neck. More precisely, decompose $X$ as

$$X = X^\circ \cup_Y N,$$

where $Y$ is unit circle bundle over $\Sigma$ with Euler number $-n$, $N$ is a tubular neighborhood of the surface $\Sigma$ (which is diffeomorphic to the disk bundle associated to $Y$), and $X^\circ$ is the complement in $X$ of the interior of $N$. Fix metrics $g_{X^\circ}$, $g_N$, and $g_Y$ for which $g_{X^\circ}$ and $g_N$ are isometric to

$$dt^2 + g_Y^2$$

in a collar neighborhood of their boundaries (where $t$ is a normal coordinate to the boundary). Let $X(T)$ denote the Riemannian manifold which is diffeomorphic to $X$ and whose metric $g_T$ is obtained from the description

$$X(T) = X^\circ \cup_{\partial X^\circ = \{-T\} \times Y} [-T, T] \times Y \cup_{\{T\} \times Y = \partial N} N;$$

i.e. $g_T|_{X^\circ} = g_{X^\circ}$, $g_T|_{[-T,T] \times Y} = dt^2 + g_Y^2$, and $g_T|_N = g_N$. For all sufficiently large $T$, there is a description of the moduli space $\mathcal{M}(X(T), \mathfrak{s}_0)$ on $X(T)$ in terms of the moduli spaces for $Y$ and the cylindrical-end, $L^2$ moduli spaces for $X^\circ$ and $N$, denoted $\mathcal{M}(X^\circ, \mathfrak{s}_0|_{X^\circ})$, and $\mathcal{M}(N, \mathfrak{s}_0|_N)$ respectively.

To understand the moduli space of $Y$, we appeal to the following result contained in [20] (we state a combination of Corollaries 5.8.5 and 5.9.1). Recall that a solution to the Seiberg-Witten equations on a four-manifold (or a three-manifold) is called *reducible* if the spinor vanishes entirely; and it is called *irreducible* otherwise. Correspondingly, we partition the moduli spaces into spaces of reducibles $\mathcal{M}^{\mathrm{red}}$ and irreducibles $\mathcal{M}^{\mathrm{irr}}$.

THEOREM 2.3 ([20]). *Let $Y$ be a circle bundle over a Riemann surface $\Sigma$ and Euler number $-n$. Let $\mathfrak{s}$ be a $\mathrm{Spin}_{\mathbb{C}}$ structure on $N$ with*

$$\langle c_1(\mathfrak{s}), [\Sigma] \rangle = k.$$

*Then, the moduli space of reducibles in $\mathfrak{s}|_Y$ is identified with the Jacobian $\mathcal{J}(\Sigma)$. This moduli space is nondegenerate unless*

$$k \equiv n \pmod{2n}.$$

*Moreover, the moduli space contains irreducibles if and only if $k$ is congruent modulo $2n$ to an integer in the range*

$$[-n-2g+2, -n-2] \cup [n+2, n+2g-2].$$

*Remark* 2.4. Strictly speaking, the above description holds for the $\nabla$-compatible moduli spaces, where $\nabla$ is a connection on $TY$ which has torsion (see [20]). Similarly, for the moduli spaces on $N$ described below, we use a connection which extends $\nabla$, which is described in Section 5. We are free to work with these moduli spaces, according to the general results of Section 4 (see Theorem 4.6).

Since $n \geq 2g$, the above result ensures that in $\mathfrak{s}_0|_Y$, the moduli space contains no irreducibles, and is diffeomorphic to the Jacobian $\mathcal{J}(\Sigma)$.

Similarly, we can completely describe the moduli space over $N$ (using a suitable connection on $TX$):

PROPOSITION 2.5. *Let $\mathfrak{s}$ be a $\mathrm{Spin}_{\mathbb{C}}$ structure, with*

$$\left|\langle c_1(\mathfrak{s}), [\Sigma] \rangle\right| < n,$$

*then for each $A \in \mathcal{M}(N, \mathfrak{s})$ which is asymptotic to a reducible, the perturbed Dirac operator $\displaystyle{\not}D_A^+$ has no kernel or cokernel. In other words, the moduli space of solutions with reducible boundary values consists entirely of reducibles, which are cut out transversally by the Seiberg-Witten equations.*

The above proposition is contained in Corollary 6.2 and Proposition 6.4.

For the $\mathrm{Spin}_{\mathbb{C}}$ structures discussed above, the moduli space of reducibles over $N$ is cut out transversely by the Seiberg-Witten equations. This is the case for no other $\mathrm{Spin}_{\mathbb{C}}$ structure: rather, the reducibles have a more subtle local Kuranishi description. For certain $\mathrm{Spin}_{\mathbb{C}}$ structures $\mathfrak{s}$, these Kuranishi models fit together to form a smooth vector bundle over $\mathcal{M}^{\mathrm{red}}(N, \mathfrak{s})$, the *obstruction bundle*. For our proof, it suffices to consider the following case:

PROPOSITION 2.6. *If*

$$\langle c_1(\mathfrak{s}), [\Sigma] \rangle = -n - 2g,$$

*then the moduli space of finite-energy monopoles consists entirely of reducibles. Moreover, for all $A \in \mathcal{M}^{\mathrm{red}}(N)$, the kernel of $\displaystyle{\not}D_A^+$ vanishes, and its cokernel*

*has complex dimension $g$. Indeed, the cokernels fit together to form a vector bundle $V$ over $\mathcal{M}^{\mathrm{red}}(N)$ whose Chern classes are given by the formula*

$$c(V) = \prod_{i=1}^{g}\Big(1 + \mu(A_i)\mu(B_i)\Big),$$

*for any standard symplectic basis $\{A_i, B_i\}_{i=1}^{g}$ for $H_1(\Sigma;\mathbb{Z})$.*

This is a combination of Proposition 5.7 and Corollary 6.3.

Now, we say what we can about the moduli space on $X^\circ$. By our assumption equation (9), Theorem 2.3 says that the moduli space of $\mathfrak{s}_0|_Y$ is smooth and consists only of reducibles. Thus, for a generic, compactly-supported two-form in $X^\circ$, the moduli space $\mathcal{M}(X^\circ)$ is smooth and contains no reducibles. Since $c_1(\mathfrak{s}_0|_Y)$ is a torsion class, the Chern-Simons-Dirac functional is real-valued. Moreover, since the $\mathcal{M}(Y,\mathfrak{s}_0|_Y)$ is connected, standard arguments show that $\mathcal{M}(X^\circ,\mathfrak{s}_0|_{X^\circ})$ is compact. (For a more detailed discussion of genericity and compactness, see [11] and [17].)

We would like to compare the Seiberg-Witten invariant of $X$ in the $\mathrm{Spin}_{\mathbb{C}}$ structure $\mathfrak{s}_0$ with the invariant in the $\mathrm{Spin}_{\mathbb{C}}$ structure

$$\mathfrak{s}_1 = \mathfrak{s}_0 - \mathrm{PD}([\Sigma]),$$

by describing both in terms of the moduli spaces over $X^\circ$. We begin with the easier case, the invariant for $\mathfrak{s}_1$. But first, we recall the geometric interpretation of the cohomology class $\mu(U)$ over the moduli space.

Given a point $x \in X$, we define a principal circle bundle over the moduli space $\mathcal{M}(X)$, the *based moduli space*, denoted $\mathcal{M}^0(X)$, consisting of the solutions to the Seiberg-Witten equations modulo gauge transformations which fix the fiber of the spinor bundle over $x$. Parallel transport allows us to place this base point anywhere on $X$ without changing the isomorphism class of the circle bundle. Let $\mathcal{L}_x \to \mathcal{M}(X^\circ)$ denote the complex line bundle associated to $\mathcal{M}^0(X^\circ)$. Then, the class $\mu(U) \in H^2(\mathcal{M}(X);\mathbb{Z})$ (used in the definition of the Seiberg-Witten invariant) is $c_1(\mathcal{L}_x)$. Over $X^\circ$, we can perform the same construction; indeed, in this case, since $\mathcal{M}(X^\circ)$ is compact, we are free to place the base point "at infinity" (i.e. we take the quotient by those gauge transformations whose limiting value at some fixed point $y \in Y$ is trivial). With these remarks in place, we turn to moduli space for $\mathfrak{s}_1$.

PROPOSITION 2.7. *The Seiberg-Witten invariant of $X$ in $\mathfrak{s}_1$ is calculated by*

$$\mathrm{SW}_{X,\mathfrak{s}_1}(a) = \langle \mu(a), [\mathcal{M}(X^\circ,\mathfrak{s}_1|_{X^\circ})]\rangle,$$

*for any $a \in \mathbb{A}(X)$.*

*Proof.* The $L^2$ theory for moduli spaces with cylindrical ends gives boundary value maps

$$\begin{aligned}\partial_{X^\circ}\colon \mathcal{M}(X^\circ, \mathfrak{s}_1|_{X^\circ}) &\longrightarrow \mathcal{M}(Y, \mathfrak{s}_1|_Y) \\ \partial_N\colon \mathcal{M}(N, \mathfrak{s}_1|_N) &\longrightarrow \mathcal{M}(Y, \mathfrak{s}_1|_Y).\end{aligned}$$

In fact, according to Proposition 2.5, $\partial_N$ is a diffeomorphism between smooth moduli spaces. Since there is no obstruction bundle over $N$, the usual gluing techniques over smooth boundary values give a fiber-product description of the moduli space $\mathcal{M}(X(T), \mathfrak{s}_1)$ for all sufficiently large $T$ (such gluing is described in [4]; see also [5] and [19]). Since $\partial_N$ is a diffeomorphism, it follows indeed that the moduli spaces $\mathcal{M}(X^\circ, \mathfrak{s}_1|_{X^\circ})$ and $\mathcal{M}(X, \mathfrak{s}_1)$ are diffeomorphic under an identification which respects the cohomology classes induced by the $\mu$-maps. By the definition of the Seiberg-Witten invariant, the proof of the proposition is complete. □

There is a similar description for the invariant in $\mathfrak{s}_0$.

PROPOSITION 2.8. *The Seiberg-Witten invariant of $X$ in $\mathfrak{s}_0$ is calculated by*

$$\mathrm{SW}_{X,\mathfrak{s}_0}(b) = \langle \mu(\xi(-\Sigma) \cdot b), [\mathcal{M}(X^\circ, \mathfrak{s}_0|_{X^\circ})] \rangle,$$

*for any $b \in \mathbb{A}(X)$.*

*Proof.* Once again, we have a fiber-product description of the moduli space on $X$. However, this time, the moduli space over $N$ is obstructed; indeed according to Proposition 2.6, we have an obstruction bundle

$$V \longrightarrow \mathcal{M}^{\mathrm{red}}(N, \mathfrak{s}_0|_N) = \mathcal{M}(N, \mathfrak{s}_0|_N).$$

In this case, gluing gives a description of $\mathcal{M}(X, \mathfrak{s}_0)$ as the zeros of a section $s$ of the bundle

$$\mathcal{M}^0(X^\circ) \times_{S^1} \partial^*_{X^\circ}(V) \longrightarrow \mathcal{M}(X^\circ).$$

Since the action of $S^1$ (viewed as constant gauge transformations) on the fibers of $V$ is the standard, weight-one action, we see that

$$\mathcal{M}^0(X^\circ) \times_{S^1} \partial^*_{X^\circ}(V) \cong \mathcal{L} \otimes \partial^*_{X^\circ}(V).$$

But the Chern class of $V$, according to Proposition 2.6, is given by

$$c(V) = \prod_{i=1}^{g}\Big(1 + \mu(A_i)\mu(B_i)\Big),$$

so

$$c(\mathcal{L} \otimes \partial^*_{X^\circ}(V)) = \prod_{i=1}^{g}\Big(1 + \mu(A_i)\mu(B_i) + \mu(U)\Big);$$

in particular,
$$e(\mathcal{L} \otimes \partial_{X^\circ}^*(V)) = \mu(\xi(-\Sigma)).$$

Thus,
$$\begin{aligned} \mathrm{SW}_{\mathfrak{s}_0}(b) &= \langle \mu(b), [s^{-1}(0)] \rangle \\ &= \langle \mu(b) \cup e(\mathcal{L} \otimes \partial_{X^\circ}^*(V)), [\mathcal{M}(X^\circ, \mathfrak{s}_0|_{X^\circ})] \rangle \\ &= \langle \mu(\xi(-\Sigma)) \cup \mu(b), [\mathcal{M}(X^\circ, \mathfrak{s}_0|_{X^\circ})] \rangle. \end{aligned}$$

This proves the proposition. □

Since the two $\mathrm{Spin}_\mathbb{C}$ structures $\mathfrak{s}_0$ and $\mathfrak{s}_1$ agree over $X^\circ$, Propositions 2.7 and 2.8 together prove Proposition 2.1. This completes the proof of the relation, save for the claims about the moduli spaces on $N$ which which we return to in Sections 4–6.

We close with some remarks on the case where $b_2^+(X) = 1$. Since the self-intersection number of $\Sigma$ is negative, while we stretch the neck along $Y$, the period point converges to a cohomology class in $\omega \in \Omega^+$ perpendicular to $\Sigma$. By choosing the perturbation two-form compactly supported in $X^\circ$, we prove the relation for all common chambers for $\mathfrak{s}$ and $\mathfrak{s} - \mathrm{PD}(\Sigma)$, which meet the line $(\omega, t)$ as $t$ varies. Since $\omega \cdot \mathrm{PD}(\Sigma) = 0$, these chambers are by definition perpendicular to $\Sigma$, and since we are free to choose $t$ arbitrarily large in magnitude, we get the relation in two different common chambers. But there are only two such chambers (see Remark 1.4); so we have proved the proposition when $b_2^+(X) = 1$.

## 3. Proof of Theorem 1.1

The symplectic Thom conjecture is a consequence of the Relation (4), and the following basic result of Taubes for the Seiberg-Witten invariants of symplectic four-manifolds.

THEOREM 3.1 (Taubes [27]). *Let $(X, \omega)$ be a closed, symplectic four-manifold. Then, for the canonical $\mathrm{Spin}_\mathbb{C}$ structure $\mathfrak{s}_0 \in \mathrm{Spin}_\mathbb{C}(X)$, we have that*

(10) $$\mathrm{SW}_{X,\mathfrak{s}_0}(1) = \pm 1;$$

*furthermore, for all other $\mathrm{Spin}_\mathbb{C}$ structures $\mathfrak{s}$, with $\mathrm{SW}_{X,\mathfrak{s}} \not\equiv 0$, we have*

(11) $$c_1(\mathfrak{s}_0) \cdot \omega \leq c_1(\mathfrak{s}) \cdot \omega.$$

*In the case where $b_2^+(X) = 1$, all Seiberg-Witten invariants should be calculated in the chamber corresponding to the perturbation $\eta = -t\omega$, for sufficiently large $t > 0$.*

*Remark* 3.2. The canonical $\mathrm{Spin}_{\mathbb{C}}$ structure $\mathfrak{s}_0$ is the one for which $W^+ \cong \Lambda^{0,0} \oplus \Lambda^{0,2}$; thus, $c_1(\mathfrak{s}_0) = -K$, where $K$ is the canonical class of the symplectic structure.

We begin with the argument when $b_2^+(X) > 1$. Suppose there is a counterexample to the theorem; i.e. there is an embedded symplectic submanifold $\Sigma \subset X$ and a homologous, smoothly-embedded submanifold $\Sigma'$ with

$$g(\Sigma') < g(\Sigma).$$

By blowing up if necessary, we can find another counterexample to the theorem for which the self-intersection number of the homology class is negative. So from now on, we assume that $-n = [\Sigma] \cdot [\Sigma] < 0$. Moreover, by attaching trivial handles to $\Sigma'$ if necessary, we can assume that

$$g(\Sigma') = g(\Sigma) - 1.$$

The adjunction formula for the symplectically embedded surface $\Sigma$ gives:

$$[\Sigma] \cdot [\Sigma] + \langle K_X, [\Sigma] \rangle = 2g(\Sigma) - 2,$$

so

$$\langle c_1(\mathfrak{s}_0), [\Sigma'] \rangle = -\langle K_X, [\Sigma] \rangle = -2g(\Sigma') - n.$$

Using equation (4) (i.e. Theorem 1.3 in the case when $g(\Sigma') > 0$ and the corresponding result of [6] in the case when $g(\Sigma') = 0$), combined with equation (10), we have that

$$\mathrm{SW}_{X,\mathfrak{s}_0-\mathrm{PD}([\Sigma'])}(\xi(-\Sigma')) = \mathrm{SW}_{X,\mathfrak{s}_0}(1) = \pm 1.$$

On the other hand,

$$c_1(\mathfrak{s}_0 - \mathrm{PD}([\Sigma'])) \cdot \omega = c_1(\mathfrak{s}_0) \cdot \omega - 2\mathrm{Vol}(\Sigma) < c_1(\mathfrak{s}_0) \cdot \omega.$$

This inequality, together with the fact that $\mathrm{SW}_{X,\mathfrak{s}_0-\mathrm{PD}(\Sigma)} \not\equiv 0$, contradicts inequality (11) of Taubes. This contradiction proves Theorem 1.1 when $b_2^+(X) > 1$.

When $b_2^+(X) = 1$, we arrange first that the self-intersection of $\Sigma$ is negative, and $g(\Sigma') = g(\Sigma) - 1$. Now choose an $\omega' \in \Omega^+(X)$ perpendicular to $\mathrm{PD}(\Sigma)$. For all large $t$, $(\omega', -t)$ is a common chamber for $\mathfrak{s}_0$ and $\mathfrak{s}_0 - \mathrm{PD}(\Sigma)$ which is orthogonal to $\mathrm{PD}(\Sigma)$ (where the relation holds) but it is also in the symplectic chamber for both $\mathrm{Spin}_{\mathbb{C}}$ structures (where Taubes' theorem holds). So the previous argument applies in this case as well.

*Remark* 3.3. Note that in the case when $b_2^+(X) > 1$ and $g(\Sigma) > 1$, the symplectic Thom conjecture also follows from our adjunction inequality together with another theorem of Taubes which states that any symplectic four-manifold with $b_2^+(X) > 1$ has Seiberg-Witten simple type; see [28].

## 4. Compactness

Typically, the Seiberg-Witten equations are viewed as equations for pairs $(A, \Phi)$, which $\Phi$ is a spinor, and $A$ is a "spin connection"; i.e. one which induces the Levi-Civita connection $\widehat{\nabla}$ on $TX$. For our calculations on the moduli spaces of $N$ (Sections 5 and 6), however, we find it convenient to use the Seiberg-Witten moduli spaces of pairs where $A$ induces not Levi-Civita on $TX$, but rather another connection (described in Section 5) which naturally arises from the bundle structure of the tubular neighborhood. The purpose of this section is to show that, in general, using alternate connections on the tangent bundle constitutes an allowable perturbation of the usual Seiberg-Witten equations, in the sense that the associated invariant is the same (see Theorem 4.6).

The crux of the matter is to derive a general Weitzenböck formula for the $\nabla$-compatible Dirac operator (for an arbitrary connection $\nabla$ on the tangent bundle), and use that to prove compactness for the corresponding moduli spaces. We begin by stating this Weitzenböck formula.

LEMMA 4.1. *Let $X$ be a four-manifold with Riemannian metric $g$, and let $\nabla$ be a compatible $\mathrm{SO}(4)$ connection on $TX$. Then, there is a vector field $\xi$ over $X$, and a pair of bundle maps*

$$\alpha \colon W^- \longrightarrow W^+;$$

$$\beta \colon W^+ \longrightarrow W^+$$

*with the property that for any $\nabla$-compatible $\mathrm{Spin}_{\mathbb{C}}$ connection $A$, we have:*

(12) $\qquad (\displaystyle{\not}D_A^+)^*(\displaystyle{\not}D_A^+) = \widehat{\nabla}_A^* \widehat{\nabla}_A + \alpha \circ \displaystyle{\not}D_A^+ - \widehat{\nabla}_{A;\xi} + \beta - \frac{1}{2}\rho_{\Lambda^+}(\mathrm{Tr}F_A^+).$

*Here, $\widehat{\nabla}_A$ denotes the covariant derivative with respect to the $\widehat{\nabla}$-compatible connection $\widehat{A}$ with the property that $\mathrm{Tr}(\widehat{A}) = \mathrm{Tr}(A)$, and $\widehat{\nabla}_{A;\xi}$ denotes covariant derivative with respect to this connection in the direction $\xi$.*

*Remark* 4.2. In the course of the proof of this lemma, we derive explicit formulas for $\alpha$, $\beta$ and $\xi$ in terms of the difference form $\widehat{\nabla} - \nabla \in \Omega^1(X, \mathfrak{so}(TX))$, but the version stated above is sufficient for our purposes.

We defer the proof of the lemma, and show how to use it to prove the following generalization of Lemma 2 from [11]:

PROPOSITION 4.3. *Let $(X, g)$ be a Riemannian four-manifold, let $\eta$ be a self-dual two-form, and let $\nabla$ be a $g$-compatible connection on $TX$. Then there is a constant $C$ depending only on $g$, $\nabla$ and $\eta$ with the property that for* *any solution* $(A, \Phi)$ *to the Seiberg-Witten equations relative to* $\nabla$ *satisfies a universal bound*

$$|\Phi|^2 \leq C.$$

*Proof.* At a point where $|\Phi|^2$ is maximal, we have that

$$\begin{aligned}
0 &\leq (d^*d - \xi)|\Phi|^2 \\
&\leq 2\langle \widehat{\nabla}_A^* \widehat{\nabla}_A \Phi, \Phi \rangle - 2\langle \widehat{\nabla}_{A;\xi} \Phi, \Phi \rangle \\
&= -\langle \beta \circ \Phi, \Phi \rangle + \langle \rho_{\Lambda^+}(i\eta)\Phi, \Phi \rangle + \langle i(\Phi \otimes \Phi^*)_\circ \Phi, \Phi \rangle \\
&= -\langle \beta \circ \Phi, \Phi \rangle + \langle \rho_{\Lambda^+}(i\eta)\Phi, \Phi \rangle - \frac{1}{2}|\Phi|^4.
\end{aligned}$$

(The first line follows from the maximum principle, the second from the Leibnitz rule, the third from the modified Weitzenböck formula (12) together with the Seiberg-Witten equations (1) and (2), and the last is a standard property of the squaring map appearing in the Seiberg-Witten equations.) This proves the proposition, for the constant

$$C = 2\max(|\beta| + |\eta|). \qquad \square$$

Thus, we have:

COROLLARY 4.4. *On a smooth, closed, oriented, Riemannian four-manifold* $X$ *equipped with a fixed connection* $\nabla$ *on* $TX$, *and self-dual two-form* $\eta$, *the moduli space of solutions to the Seiberg-Witten equations compatible with* $\nabla$ *and perturbed by* $\eta$ *is compact.*

The proof of Corollary 4.4 is standard, given the bound from Proposition 4.3, together with elliptic regularity (see [11] for a similar statement).

Another consequence of the Weitzenböck formula is the following fundamental unique continuation principle.

PROPOSITION 4.5. *For* $X$, $g$, *and* $\nabla$ *as in Proposition* 4.3, *let* $(A, \Phi)$ *be a solution to the* $\eta$-*perturbed,* $\nabla$-*compatible Seiberg-Witten equations. Then, if* $\Phi$ *vanishes identically on any open subset of* $X$, *then indeed* $\Phi$ *must vanish identically over* $X$.

*Proof.* If $\Phi$ satisfies $\displaystyle{\not}D_A^+ \Phi \equiv 0$, then the Weitzenböck formula says that

$$\widehat{\nabla}_A^* \widehat{\nabla}_A \Phi = \widehat{\nabla}_{A;\xi} \Phi + \beta(\Phi) - \frac{1}{2}\rho_{\Lambda^+}(\mathrm{Tr}F_A^+)\Phi;$$

in particular, $\Phi$ satisfies an inequality of the form

$$|\widehat{\nabla}_A^* \widehat{\nabla}_A \Phi|^2 \leq M(|\Phi|^2 + |\widehat{\nabla}_A \Phi|^2)$$

for some constant $M$. The result then follows from a general result of Aronszajn (see [1]). $\qquad \square$

Given the compactness and unique continuation results established above (Corollary 4.4 and Proposition 4.5 respectively), the usual arguments for proving metric independence of the Seiberg-Witten invariants using the Levi-Civita connection (see for example [11] or [16]) apply *mutas mutandis* to prove that the invariants using a metric $g$ and any $g$-compatible $SO(4)$ connection $\nabla$ over $TX$ are actually independent of the pair $(g, \nabla)$, save for a chamber dependence when $b_2^+(X) = 1$. Explicitly, letting $\mathcal{M}(X, g, \nabla, \eta, \mathfrak{s})$ denote the moduli space of gauge equivalence classes of pairs $(A, \Phi)$ satisfying the $\eta$-perturbed Seiberg-Witten equations, where now $A$ is a connection inducing $\nabla$ on $TX$, we have the following theorem:

THEOREM 4.6. *Let $X$ be a smooth, oriented, closed four-manifold. For any metric $g$, any $g$-compatible $SO(4)$-connection $\nabla$ on $TX$, there is a dense set of two-forms $\eta \in \Omega^+(X, \mathbb{R})$ with the property that the $\eta$-perturbed, $\nabla$-compatible Seiberg-Witten moduli space $\mathcal{M}(X, g, \nabla, \eta, \mathfrak{s})$ is a compact, smooth moduli space of expected dimension $d(\mathfrak{s})$. Moreover, after orienting the moduli space with a homology orientation, the homological pairing*

$$\langle \mu(a), [\mathcal{M}(X, g, \nabla, \eta, \mathfrak{s})]\rangle$$

*is independent of $\nabla$. Indeed, if $b_2^+(X) > 1$, it is independent of $g$ and $\eta$ as well. If $b_2^+(X) = 1$, it depends on $g$ and $\eta$ only through the chamber of the associated period point.*

*Proof.* Since most of this statement is standard, we content ourselves with an outline.

For a metric $g$ and $g$-compatible connection $\nabla$ on $TX$, let $\mathfrak{M}^{\mathrm{irr}}(X, g, \nabla, \mathfrak{s})$ denote the universal irreducible moduli space, the space of gauge equivalence classes of triples $[A, \Phi, \eta]$, where $A$ is a $\nabla$-compatible spinor connection, $\Phi$ is a spinor which does not vanish identically, and the pair $(A, \Phi)$ satisfies the $\eta$-perturbed equations. One first shows that for any $(g, \nabla)$ as above, $\mathfrak{M}^{\mathrm{irr}}(X, g, \nabla, \mathfrak{s})$, is a smooth Hilbert manifold. This follows by applying the arguments of Section 2 of [11], bearing in mind that unique continuation (Proposition 4.5) still holds. Thus, Sard-Smale theory ([23], see also [4]) shows that for a fixed pair $(g, \nabla)$, there is a dense (Baire) set of $\eta$ for which the associated moduli space $\mathcal{M}^{\mathrm{irr}}(X, g, \nabla, \eta, \mathfrak{s})$ is smooth. Similarly, the Sard-Smale transversality theorem shows that given any two triples $(g_0, \nabla_0, \eta_0)$ and $(g_1, \nabla_1, \eta_1)$ whose associated moduli spaces are smooth, and any path connecting those triples, there is another path $(g_t, \nabla_t, \eta_t)$ for $t \in [0, 1]$ (which we can take to be arbitrarily close to the original path) with the property that the corresponding one-parameter family of moduli spaces forms a smooth cobordism between $\mathcal{M}^{\mathrm{irr}}(X, g_0, \nabla_0, \eta_0, \mathfrak{s})$ and $\mathcal{M}^{\mathrm{irr}}(X, g_1, \nabla_1, \eta_1, \mathfrak{s})$. This cobordism is compact provided that there are no reducible solutions in the moduli spaces spaces $\mathcal{M}(X, g_t, \nabla_t, \eta_t, \mathfrak{s})$ for $t \in [0, 1]$. As usual, reducibles occur in a codimension $b_2^+(X)$ subspace in the product space of metrics and forms, so they can be avoided when $b_2^+(X) > 1$, and when $b_2^+(X) = 1$, they occur precisely when the period points associated to $(g_t, \eta_t)$ cross walls. (In particular, the condition that a moduli space contains reducibles makes no reference to the connection $\nabla$ on $TX$; it is the condition that some line bundle – the determinant of $W^+$ – has a curvature form representative with specified self-dual part.)

Given this smooth cobordism, the statement about invariants follows in the usual manner. $\square$

Our remaining goal in this section is to prove Lemma 4.1. The derivation will employ the following standard fact about the Dirac operator (coupled to the Levi-Civita connection).

LEMMA 4.7. *Let $\widehat{\nabla}_A$ be a connection on the spinor bundle compatible with the Levi-Civita connection $\widehat{\nabla}$ on $TX$, and let $\theta$ be a smooth one-form. Then, the anti-commutator*

$$\{\widehat{\slashed{D}}_A, \rho(\theta)\} = \widehat{\slashed{D}}_A \circ \rho(\theta) + \rho(\theta) \circ \widehat{\slashed{D}}_A$$

*satisfies the following relation:*

(13) $$\{\widehat{\slashed{D}}_A, \rho(\theta)\} = -2\widehat{\nabla}_{A;\theta^\flat} + \rho((d+d^*)\theta),$$

*where $\theta^\flat$ is the vector field which is dual to $\theta$.*

*Remark* 4.8. Note that the Clifford action $\rho$ induces an action of the entire exterior algebra $\sum_{p=1}^4 \Lambda^p T^* X$ on the full spinor bundle $W = W^+ \oplus W^-$. We denote this extended action by $\rho$ as well (as in the right-hand side of equation (13)). Occasionally, we will write $\rho_{\Lambda^p}$ for the restriction of this action to $\Lambda^p T^* X$ in the interest of clarity.

A proof of the above can be found in [3, p. 122]. We now derive equation (12).

*Proof of Lemma* 4.1. We prove the Weitzenböck formula for a connection $A$ coupled to $\nabla$, using the usual Weitzenböck formula for connections coupled to the Levi-Civita connection $\widehat{\nabla}$.

First, observe that there are forms $\mu \in \Omega^1(X)$, $\nu \in \Omega^3(X)$ with the property that

(14) $$\widehat{\slashed{D}}_A = \slashed{D}_A + \rho(\mu) + \rho(\nu).$$

Note that this equation is for the Dirac operator $\slashed{D}_A = \slashed{D}_A^+ \oplus \slashed{D}_A^-$ acting on the full spinor bundle $W = W^+ \oplus W^-$. The forms $\mu$ and $\nu$ are extracted from the connection one-form[1] $\omega = \widehat{\nabla} - \nabla \in \Omega^1(X, \mathfrak{so}(TX))$ as follows. Recall that the natural map

$$i \colon \mathfrak{so}(TX) \longrightarrow \Lambda^2(T^*X)$$

defined by

$$i(a_j^k) = \frac{1}{4} \sum_{j,k} a_j^k \theta^j \wedge \theta^k$$

is a vector space isomorphism, and that the action of $\mathfrak{so}(TX)$ on the Clifford bundle $W$ is modeled on $\rho_{\Lambda^2} \circ i$ (see [9]). Thus, via this isomorphism, we can view the difference one-form as an element $\omega \in \Gamma(X, (\Lambda^1 \otimes \Lambda^2)T^*X)$, and the connection $\widehat{\nabla}_A - \nabla_A$ is the one-form induced by Clifford multiplying the $\Lambda^2$ component of $\omega$. Thus,

$$\widehat{\slashed{D}}_A - \slashed{D}_A = \rho_{\Lambda^1 \otimes \Lambda^2}(\omega),$$

where

$$\rho_{\Lambda^1 \otimes \Lambda^2} \colon \Lambda^1 \otimes \Lambda^2 \longrightarrow \mathrm{End}(W)$$

denotes the linear map with the property that

$$\rho_{\Lambda^1 \otimes \Lambda^2}(\theta \otimes \gamma) = \rho_{\Lambda^1}(\theta) \circ \rho_{\Lambda^2}(\gamma)$$

for each $\theta \in \Lambda^1$, $\gamma \in \Lambda^2$. Since moreover

$$\rho_{\Lambda^1 \otimes \Lambda^2}(\theta \otimes \gamma) = \rho_{\Lambda^1}(\iota_{\theta^\flat} \gamma) + \rho_{\Lambda^3}(\theta \wedge \gamma),$$

where $\iota_{\theta^\flat}$ denotes contraction, we can find the forms $\mu$ and $\nu$ appearing in equation (14).

Restricting attention to $W^+$, the induced action of any three-form $\nu$ agrees with the action of the one-form $*\nu$; so it follows from equation (14), that

$$\widehat{\slashed{D}}_A^+ = \slashed{D}_A^+ + \rho(\gamma),$$

where $\gamma$ is the the one-form $\gamma = \mu + *\nu$. Thus, with the help of the Anti-Commutator Formula (13) and the fact that $\widehat{\slashed{D}}_A$ is self-adjoint, we have for for any one-form $\theta$ that

$$\begin{aligned}(\slashed{D}_A^+)^* \circ \rho(\theta) &= ((\widehat{\slashed{D}}_A^+)^* + \rho(\gamma)) \circ \rho(\theta) \\ &= -\rho(\theta) \circ \widehat{\slashed{D}}_A^+ - 2\widehat{\nabla}_{A;\theta^\flat} + \rho((d^* + d)\theta) + \rho(\gamma)\rho(\theta) \\ &= -\rho(\theta) \circ \slashed{D}_A^+ - 2\widehat{\nabla}_{A;\theta^\flat} + [\rho(\gamma), \rho(\theta)] + \rho((d^* + d)\theta);\end{aligned}$$

in particular,

$$(\slashed{D}_A^+)^* \circ \rho(\gamma) = -\rho(\gamma) \circ \slashed{D}_A^+ - 2\widehat{\nabla}_{A;\gamma^\flat} + \rho((d^* + d)\gamma).$$

---

[1] In keeping with standard notation for the Cartan formalism (see for instance [24]), we let $\omega$ denote the connection form for connections on $TX$. This should not be confused with the symplectic form from Section 3. Indeed, we will not make use of any symplectic forms for the rest of this paper.

Thus,

$$\begin{aligned}(\widehat{\slashed{D}}_A^+)^* \circ (\widehat{\slashed{D}}_A^+) &= ((\slashed{D}_A^+)^* - \rho(\gamma)) \circ (\slashed{D}_A^+ + \rho(\gamma)) \\ &= (\slashed{D}_A^+)^* \circ \slashed{D}_A^+ - \rho(\gamma) \circ \slashed{D}_A^+ + (\slashed{D}_A^+)^* \circ \rho(\gamma) - \rho(\gamma) \circ \rho(\gamma) \\ &= (\slashed{D}_A^+)^* \circ \slashed{D}_A^+ - 2\rho(\gamma) \circ \slashed{D}_A^+ - 2\widehat{\nabla}_{A;\gamma^\flat} + \rho((d^* + d)\gamma) + |\gamma|^2.\end{aligned}$$

Substituting in the usual Weitzenböck formula (see [9])

$$(\widehat{\slashed{D}}_A^+)^* \widehat{\slashed{D}}_A^+ \Phi = \widehat{\nabla}_A^* \widehat{\nabla}_A \Phi + \frac{s}{4}\Phi - \frac{1}{2}\rho_{\Lambda^+}(\mathrm{Tr} F_A^+)\Phi,$$

(where $s$ is the scalar curvature function of $X$) and rearranging terms, we obtain a formula of the shape given in equation (12). □

## 5. A connection on $N$

In this section, we describe a felicitous connection $\nabla$ on the tangent bundle $TN$ and show that its corresponding Dirac operator is $\sqrt{2}(\overline{\partial}_A + \overline{\partial}_A^*)$. With this holomorphic interpretation in hand, we begin to analyze the moduli spaces over $N$ encountered in Section 2, an enterprise which we complete in Section 6. Our constructions in this section naturally extend the analogous ones from [20].

We begin by constructing an appropriate metric on $N$. First, choose a connection on $Y$, the circle bundle over $\Sigma$, and denote the connection one-form by $\varphi$. Choose a smooth, real-valued function

$$f: (0, \infty) \longrightarrow (0, \infty),$$

with

$$f(t) = \begin{cases} f(t) = t & \text{for } t \in (0, \tfrac{1}{2}] \\ f(t) = 1 & \text{for } t \in [1, \infty). \end{cases}$$

Using this function, define a metric on the punctured disk $D - \{0\} \cong (0, \infty) \times S^1$ (with polar coordinates $(t, \theta)$) given by

$$g = dt^2 + f(t)^2 d\theta^2.$$

Obviously, this metric extends across the origin in $D$ (the metric on $(0, \tfrac{1}{2}) \times S^1$ is the flat metric on the punctured disk). Denote the induced metric by $g_D$. There is a corresponding metric on $N - \Sigma \cong (0, \infty) \times Y$

$$g_N = \pi^* g_\Sigma + dt^2 + f(t)^2 \varphi^2.$$

Once again, this metric extends smoothly over all of $N$.

The connection on $Y$ induces a splitting at each $x \in N$:

$$T_x N \cong \pi^*(T_{\pi(x)}\Sigma) \oplus T_x \pi^{-1}(x).$$

(There is a corresponding splitting of the cotangent bundle, where the first factor is generated by forms which pull up from $\Sigma$, and the second factor is generated by $dt$ and $\varphi$.) Letting $\nabla_\Sigma$ and $\nabla_D$ denote the Levi-Civita connections on $\Sigma$ and $D$ respectively, we consider the $\mathfrak{so}(4)$ connection $\nabla$ on $TN$ given by

$$\nabla = \nabla_\Sigma \oplus \nabla_D$$

with respect to the above splitting. Clearly, over the cylindrical region $[1, \infty) \times Y$, the connection agrees with the cylindrical connection in [20]; in particular, it has torsion. We would like to compare it with the Levi-Civita connection, but first we introduce some notions.

*Definition* 5.1. A form $\omega \in \Gamma(X, \Lambda^1 \otimes \Lambda^2)$ is called *completely off-diagonal* if there is a three-form $\nu$ so that

$$\rho_{\Lambda^1 \otimes \Lambda^2}(\omega) = \rho_{\Lambda^3}(\nu),$$

and for each one-form $\theta \in \Omega^1(X)$,

(15) $$2\rho_{\Lambda^2}(\iota_{\theta^\flat}\omega) + \{\rho_{\Lambda^1 \otimes \Lambda^2}(\omega), \rho(\theta)\} = 0.$$

In the above expression, the contraction is performed on the $\Lambda^1$ factor of $\omega$; i.e.

$$\iota_\xi(\theta \otimes \gamma) = (\iota_\xi\theta)\gamma,$$

for $\theta \in \Omega^1$, $\gamma \in \Omega^2$, and $\xi \in \Gamma(X, TX)$.

The definition is perhaps justified by the following:

LEMMA 5.2. *A form $\omega \in \Gamma(X, \Lambda^1 \otimes \Lambda^2)$ is completely off-diagonal if it lies in the subbundle generated by*

$$\theta^1 \otimes \theta^2 \wedge \theta^3 - \theta^3 \otimes \theta^1 \wedge \theta^2 + \theta^2 \otimes \theta^3 \wedge \theta^1,$$

*where $\{\theta^1, \theta^2, \theta^3\}$ are local orthonormal one-forms.*

*Proof.* Let $e_1, e_2, e_3, e_4$ be an orthonormal frame, and let $\theta^i$ be the dual of $e_i$. To check that

$$\omega = \theta^1 \otimes \theta^2 \wedge \theta^3 - \theta^3 \otimes \theta^1 \wedge \theta^2 + \theta^2 \otimes \theta^3 \wedge \theta^1$$

is completely off-diagonal, it suffices to verify equation (15) for the one-forms $\theta = \theta^i$ for $i = 1, \ldots, 4$. This is a straightforward computation. $\square$

LEMMA 5.3. *Let $\nabla$ be the connection on $TN$ described above, let $\widehat{\nabla}$ be the Levi-Civita connection on $TN$, and let $\Theta$ be the $(1,1)$ form induced by the metric and complex structure on $N$. Then,*

$$\nabla\Theta = 0,$$

SYMPLECTIC THOM CONJECTURE 113*and the difference form*

$$\omega = \widehat{\nabla} - \nabla$$

*is completely off-diagonal.*

*Remark* 5.4. In the above statement, and indeed throughout the rest of this section, we are implicitly using the identification

$$\Omega^1(\mathfrak{so}(TN)) \cong \Gamma(N; \Lambda^1 \otimes \Lambda^2)$$

induced by the isomorphism $i$ from the proof of Lemma 4.1.

*Proof.* The fact that $\Theta$ is covariantly constant follows immediately from the definition of $\nabla$, and the fact that $g_\Sigma$ and $g_D$ are Kähler metrics. To verify the off-diagonality, we employ Cartan's method of moving frames (see [24]). Let $\theta^1, \theta^2$ be a local orthonormal coframe on $T^*\Sigma$, and let $\overline{\omega}^i_j$ ($i, j = 1, 2$) denote its connection one-forms. Then, $\theta^1, \theta^2$ (viewed as a partial coframe on $TN$) can be completed to an orthonormal coframe by adjoining $\theta^3 = dt$ and $\theta^4 = f(t)\varphi$. Note that

$$(16) \qquad d\theta^4 = (\frac{\partial}{\partial t}f)f^{-1}\theta^3 \wedge \theta^4 + f\kappa \theta^1 \wedge \theta^2,$$

where $\kappa$ is a function on $\Sigma$ (which measures the curvature of the connection on $Y$). The Levi-Civita connection matrix $\widehat{\omega}$ reads

$$\begin{pmatrix} \overline{\omega}^1_1 & \overline{\omega}^1_2 - \frac{1}{2}\kappa f\theta^4 & 0 & -\frac{1}{2}\kappa f\theta^2 \\ \overline{\omega}^2_1 + \frac{1}{2}f\kappa\theta^4 & \overline{\omega}^2_2 & 0 & \frac{1}{2}f\kappa\theta^1 \\ 0 & 0 & 0 & -\frac{\partial}{\partial t}(\log f)\theta^4 \\ \frac{1}{2}\kappa f\theta^2 & -\frac{1}{2}\kappa f\theta^1 & \frac{\partial}{\partial t}(\log f)\theta^4 & 0 \end{pmatrix}.$$

(It is straightforward to verify this assertion given equation (16); one need only check that

$$d\theta = -\widehat{\omega} \wedge \theta.)$$

With respect to this coframe, the connection matrix for $\nabla$ is

$$\begin{pmatrix} \overline{\omega}^1_1 & \overline{\omega}^1_2 & 0 & 0 \\ \overline{\omega}^2_1 & \overline{\omega}^2_2 & 0 & 0 \\ 0 & 0 & 0 & -\frac{\partial}{\partial t}(\log f)\theta^4 \\ 0 & 0 & \frac{\partial}{\partial t}(\log f)\theta^4 & 0 \end{pmatrix}.$$

Thus, we can write the difference form $\omega = \widehat{\nabla} - \nabla$ as

$$\omega = \frac{\kappa f}{4}(\theta^4 \otimes \theta^1 \wedge \theta^2 + \theta^2 \otimes \theta^1 \wedge \theta^4 - \theta^1 \otimes \theta^2 \wedge \theta^4),$$

which is completely off-diagonal according to Lemma 5.2. $\square$

In light of the above calculation, the connection $\nabla$ on $TN$ behaves morally much like the the Levi-Civita connection on a Kähler manifold. We make this more precise in the following lemma, whose proof should be compared to Taubes' description of the Dirac operator on a symplectic manifold (see [27]).

LEMMA 5.5. *Let $X$ be a complex surface with metric $g$ and associated $(1,1)$ form $\Theta$. Let $\nabla$ be an $\mathfrak{so}(4)$ connection on $TX$ and suppose that*

1. $\nabla \Theta = 0$

2. *and the difference form $\widehat{\nabla} - \nabla = \omega$ is completely off-diagonal.*

*Then, for any $\nabla$-compatible connection $A$ on the spinor bundle*
$$\Lambda^{0,0}E \oplus \Lambda^{0,2}E = E \oplus K_X^{-1}(E),$$
*the Dirac operator $\slashed{D}_A^+$ is identified with the operator $\sqrt{2}(\overline{\partial}_A + \overline{\partial}_A^*)$, where the partial connection $\overline{\partial}_A$ is induced by restricting $A$ to the line bundle $E$.*

*Proof.* Consider the canonical $\mathrm{Spin}_{\mathbb{C}}$ structure whose spinor bundles are
$$W_0^+ = \Lambda^{0,0} \oplus \Lambda^{0,2} = \mathbb{C} \oplus K_X^{-1} \quad \text{and} \quad W_0^- = \Lambda^{0,1},$$
endowed with the Clifford action given by $\sqrt{2}$ times the symbol of $(\overline{\partial} + \overline{\partial}^*)$. Let $\Phi_0$ be a nonvanishing section of the $\mathbb{C}$ factor. The $\mathrm{Spin}_{\mathbb{C}}$ structure admits a $\nabla$-spinorial connection $A_0$ characterized by the property that the $\Lambda^{0,0}$ part of $\nabla_{A_0}\Phi_0$ vanishes; i.e.
$$\nabla_{A_0}\Phi_0 \in \Omega^1(X; \Lambda^{0,2}) \subset \Omega^1(X; W_0^+).$$
We prove the lemma first for the connection $A = A_0$.

Since $\Theta$ is covariantly constant, it follows that any $\nabla$-spinorial connection connection must preserve the splitting of $W_0^+$ (which can be described as the $\pm 2i$ eigenspaces of Clifford multiplication by $\Theta$). Thus, it follows that

(17) $$\nabla_{A_0}\Phi_0 \equiv 0, \quad \text{so} \quad \slashed{D}_{A_0}^+ \Phi_0 \equiv 0,$$

and indeed
$$\slashed{D}_{A_0}^+(f\Phi_0) = df \cdot \Phi_0 + f\slashed{D}_{A_0}^+ \Phi_0 = df \cdot \Phi_0 = \sqrt{2}(\overline{\partial}f)\Phi_0.$$

This gives the identification of $\slashed{D}_{A_0}^+|_{\Lambda^{0,0}}$ with $\sqrt{2}\,\overline{\partial}\colon \Omega^{0,0} \longrightarrow \Omega^{0,1}$. It remains to show that
$$\slashed{D}_{A_0}^+|_{\Lambda^{0,2}} = \sqrt{2}\,\overline{\partial}^*\colon \Omega^{0,2} \to \Omega^{0,1}$$

or, equivalently, that

(18) $$\Pi_{\Lambda^{0,2}} \slashed{D}_{A_0}^- = \sqrt{2}\,\overline{\partial}\colon \Omega^{0,1} \to \Omega^{0,2}.$$

(Here we are implicitly using the fact that $\slashed{D}_{A_0}$ is self-adjoint, which in turn follows from the fact that $\widehat{\slashed{D}}_{A_0}$ is, and so is the difference $\widehat{\slashed{D}}_{A_0} - \slashed{D}_{A_0}$, as the latter can be written as Clifford multiplication by a three-form by the hypothesis on $\omega$.) Moreover, on a complex manifold, of all first-order operators from $\Omega^{0,1} \to \Omega^{0,2}$ with the same symbol as the $\overline{\partial}$-operator, the $\overline{\partial}$ operator itself is characterized by the property that that it vanishes on the image of $\overline{\partial}: \Omega^{0,0} \longrightarrow \Omega^{0,1}$ (this can be seen, for instance, by writing out forms in local coordinates). Now, it is easy to see that the symbols of the operators compared in equation (18) agree, so the equation is established once we verify that the restriction of $\Pi_{\Lambda^{0,2}}\slashed{D}^-_{A_0}$ to the image of $\slashed{D}^+_{A_0}|_{\Lambda^{0,0}}$ vanishes. To see this, observe first that

$$\begin{aligned}
\{\slashed{D}_{A_0}, \rho(df)\} &= \{\widehat{\slashed{D}} - \rho_{\Lambda^1 \otimes \Lambda^2}(\omega), \rho(df)\} \\
&= \{\widehat{\slashed{D}}, \rho(df)\} - \{\rho_{\Lambda^1 \otimes \Lambda^2}(\omega), \rho(df)\} \\
&= (d^* + d)(df) - 2\widehat{\nabla}_{df^\flat} - \{\rho_{\Lambda^1 \otimes \Lambda^2}(\omega), \rho(df)\} \\
&= d^*df - 2\nabla_{df^\flat} - 2\rho(\iota_{df^\flat}\omega) - \{\rho_{\Lambda^1 \otimes \Lambda^2}(\omega), \rho(df)\},
\end{aligned}$$

Projecting onto $\Lambda^{0,2}$, we lose the first term; and applying the above to $\Phi_0$ (bearing in mind equation (17)), we get that

$$\Pi_{\Lambda^{0,2}} \slashed{D}^-_{A_0} \slashed{D}^+_{A_0}(f\Phi_0) = -\Pi_{\Lambda^{0,2}}(2\rho(\iota_{df^\flat}\omega) + \{\rho_{\Lambda^1 \otimes \Lambda^2}(\omega), \rho(df)\}).$$

But the latter vanishes by the hypothesis that $\omega$ is completely off-diagonal, completing the verification of equation (18).

We have shown that for the canonical $\text{Spin}_\mathbb{C}$ structure, the Dirac operator $\slashed{D}^+_{A_0}$ is identified with $\sqrt{2}(\overline{\partial} + \overline{\partial}^*)$. The corresponding statement for all other Dirac operators coupled to $\nabla$ follows immediately (by taking tensor products). □

Put together, the lemmas of this section amount to the following.

PROPOSITION 5.6. *The Dirac operator over $N$ coupled to $\nabla$ for the $\text{Spin}_\mathbb{C}$ structure $E \oplus K_N^{-1}E$ takes the form*

$$\slashed{D}^+_A = \sqrt{2}(\overline{\partial}_A + \overline{\partial}^*_A).$$

In the usual manner, we see then that the Seiberg-Witten equations over $N$ read, just as in the Kähler case [29]:

$$\begin{aligned}
\Lambda \text{Tr} F_A &= \frac{i}{2}(|\alpha|^2 - |\beta|^2) \\
\text{Tr} F_A^{0,2} &= \overline{\alpha} \otimes \beta \\
\overline{\partial}_A \alpha + \overline{\partial}^*_A \beta &= 0.
\end{aligned}$$

For finite energy solutions, decay estimates justify the usual integration-by-parts which shows that one of $\alpha$ or $\beta$ must vanish (see [20], where this is done on a cylinder). Hence, we can conclude the following:

PROPOSITION 5.7. *For any* $\mathrm{Spin}_{\mathbb{C}}$ *structure* $\mathfrak{s}$ *over* $N$ *with*

$$\left|\langle c_1(\mathfrak{s}), [\Sigma]\rangle\right| \geq n + 2g,$$

*the moduli space of finite energy solutions* $\mathcal{M}(N, \mathfrak{s})$ *consists entirely of reducibles.*

*Proof.* Write the $\mathrm{Spin}_{\mathbb{C}}$ structure as $E \oplus \Lambda^{0,2}E$. Suppose that

$$\langle c_1(\mathfrak{s}), [\Sigma]\rangle \leq -n - 2g.$$

Then, it follows that $\beta$ must vanish, since

$$\begin{aligned}\langle c_1(\mathfrak{s}), [\Sigma]\rangle &= \frac{i}{2\pi}\int_\Sigma \mathrm{Tr} F_A|_\Sigma \\ &= \frac{i}{2\pi}\int_\Sigma \Lambda\mathrm{Tr} F_A(*_\Sigma 1) \\ &= \tfrac{-1}{4\pi}\int_\Sigma |\alpha|^2 \text{ or } \tfrac{1}{4\pi}\int_\Sigma |\beta|^2,\end{aligned}$$

which we assumed to be negative. Thus, $\alpha|_\Sigma$ is a holomorphic section of $E$. But our assumption guarantees that

$$\deg(E) = \frac{1}{2}(\langle c_1(\mathfrak{s}), [\Sigma]\rangle + n + 2g - 2) < 0,$$

a bundle which admits no holomorphic sections; so it is necessary then that $\alpha = 0$ as well. The case where

$$\langle c_1(\mathfrak{s}), [\Sigma]\rangle \geq n + 2g$$

follows in an analogous fashion. □

## 6. The obstruction bundle over a neighborhood

Let $\Sigma$ be a surface of genus $g$, and $N$ denote the four-manifold which is the total space of a complex line bundle $L$ over $\Sigma$ with first Chern number $-n < 0$ (given a cylindrical-end metric $g_N$ as in §5). The goal of this section is to analyze kernel and the cokernel of the modified $\mathrm{Spin}_{\mathbb{C}}$-Dirac operator coupled to an integrable connection $A$, to complete the arguments outlined in Section 2.

A $\mathrm{Spin}_{\mathbb{C}}$ structure on $N$ is specified by giving a Hermitian line bundle $E$ over $N$, with the convention that the bundle of spinors are given by:

$$W^+ = E \oplus K_N^{-1}(E) \qquad W^- = \Lambda^{0,1}(E).$$

Numerically, a $\text{Spin}_{\mathbb{C}}$ structure is specified by the first Chern number of its determinate line bundle
$$k = \langle c_1(W^+), [\Sigma] \rangle,$$
or, equivalently, by the first Chern number of $E$,
$$e = \langle c_1(E), [\Sigma] \rangle.$$
By the adjunction formula, these quantities are related by the formula:
$$k = 2e - 2g + 2 - n.$$

In the following, let $E_0 = E|_\Sigma$. Throughout this discussion, $A$ will be an integrable connection which is asymptotic on the end to a reducible solution to the Seiberg-Witten equations over $Y$ (e.g. $A$ could be a finite-energy reducible over $N$). The following proposition gives a calculation of the $L^2$ kernel and cokernel of the modified Dirac operator
$$\sqrt{2}(\overline{\partial}_A + \overline{\partial}_A^*) \colon W^+ \longrightarrow W^-$$
(see Proposition 5.6).

PROPOSITION 6.1. *Let $\ell = \lfloor -\frac{k}{2n} - \frac{1}{2} \rfloor$, the greatest integer smaller than $-\frac{k}{2n} - \frac{1}{2}$. Assume the boundary value of $A$ is smooth; e.g. assume $n$ does not divide $k$. Then, if $\ell \geq 0$,*
$$\mathrm{Ker}(\overline{\partial}_A + \overline{\partial}_A^*) \cap L^2(W^+) \cong \sum_{j=0}^{\ell} H^0(\Sigma, E_0 \otimes L^{\otimes j}),$$
$$\mathrm{Coker}(\overline{\partial}_A + \overline{\partial}_A^*) \cap L^2(W^-) \cong \sum_{j=0}^{\ell} H^1(\Sigma, E_0 \otimes L^{\otimes j}),$$
*and if $\ell < 0$,*
$$\mathrm{Ker}(\overline{\partial}_A + \overline{\partial}_A^*) \cap L^2(W^+) \cong \sum_{j=1}^{-\ell-1} H^1(\Sigma, E_0 \otimes L^{\otimes -j}),$$
$$\mathrm{Coker}(\overline{\partial}_A + \overline{\partial}_A^*) \cap L^2(W^-) \cong \sum_{j=1}^{-\ell-1} H^0(\Sigma, E_0 \otimes L^{\otimes -j}).$$

*Proof.* Since the boundary value of $A$ is smooth, we can "unroll" the Dirac operator, giving a holomorphic interpretation of if its kernel (resp. cokernel) as the even (resp. odd) cohomology groups of the complex
$$\Omega^{0,0}(E) \xrightarrow{\overline{\partial}_A} \Omega^{0,1}(E) \xrightarrow{\overline{\partial}_A} \Omega^{0,2}(E).$$
These cohomology groups can be analyzed in the manner of Atiyah-Patodi-Singer [2] (see also [20, §§9 and 10], where this is done in great detail), to obtain an identification with the cohomology groups of the associated ruled surface $R = \mathbb{P}(\mathbb{C} \oplus L)$ with values in a line bundle $\widehat{E}$, which is specified by

$$\widehat{E}|_{\Sigma_-} \cong E_0; \qquad \widehat{E}|_{\Sigma_+} \cong \lfloor \tfrac{K}{2} \rfloor;$$

where $\Sigma_\pm$ are the two curves at infinity in the ruled surface, with

$$\Sigma_\pm \cdot \Sigma_\pm = \pm n,$$

and $\lfloor \tfrac{K}{2} \rfloor$ is the line bundle over $\Sigma$ characterized by the two properties that

$$g - 1 \geq \deg(\lfloor \tfrac{K}{2} \rfloor) > g - 1 - n$$

and

$$\deg(\lfloor \tfrac{K}{2} \rfloor) \equiv \deg(E) \pmod{n}.$$

Thus, if $F$ is the fiber in the ruling, then, as

$$\mathrm{PD}(F) = -\frac{\mathrm{PD}(\Sigma_-)}{n} + \frac{\mathrm{PD}(\Sigma_+)}{n},$$

we have that

$$\begin{aligned}\langle c_1(\widehat{E}), [F]\rangle &= \left(-\frac{e}{n} + \frac{\deg(\lfloor \tfrac{K}{2} \rfloor)}{n}\right) \\ &= \lfloor \frac{-e + g - 1}{n} \rfloor.\end{aligned}$$

Thus, the cohomology groups over $N$ of $E$ are identified with the cohomology groups over $R$ of

$$\widehat{E} = \pi^*(E_0) \otimes \mathcal{O}(\ell),$$

where $\mathcal{O}(\ell)$ denotes (fiber-wise) $\ell^{\mathrm{th}}$ twisting sheaf. We sketch the calculation, and cite [8] for details.

The Leray spectral sequence gives

$$H^i(\Sigma, \mathrm{R}^j \pi_*(\widehat{E})) \Rightarrow H^{i+j}(R, \widehat{E}).$$

By the "projection formula,"

$$\mathrm{R}^j \pi_*(\pi^*(E_0) \otimes \mathcal{O}(\ell)) = E_0 \otimes \mathrm{R}^j \pi_*(\mathcal{O}(\ell)).$$

If $\ell \geq 0$, then

$$\begin{aligned}\mathrm{R}^0 \pi_*(\mathcal{O}(\ell)) &\cong \mathrm{Sym}^\ell(\mathbb{C} \oplus L) = \sum_{j=0}^{\ell} L^j; \\ \mathrm{R}^1 \pi_*(\mathcal{O}(\ell)) &= 0.\end{aligned}$$

(Here, $\mathrm{Sym}^j \mathcal{F}$ denotes that $j$-fold symmetric product of the sheaf $\mathcal{F}$.) On the other hand, if $\ell < 0$, then

$$\mathrm{R}^0 \pi_*(\mathcal{O}(\ell)) = 0,$$

and the relative version of Serre duality gives

$$\begin{aligned} R^1\pi_*(\mathcal{O}(\ell)) &\cong \pi_*(\mathcal{O}(-\ell-2))^* \otimes (\Lambda^2(\mathbb{C} \oplus L))^* \\ &= \operatorname{Sym}^{-\ell-2}(\mathbb{C} \oplus L)^* \otimes L^* \\ &= \big(\sum_{j=0}^{-\ell-2} L^{-j}\big) \otimes L^{-1} \\ &= \sum_{j=1}^{-\ell-1} L^{-j}. \end{aligned}$$

(Here, $\mathcal{F}^*$ denotes the dual sheaf of $\mathcal{F}$.) The proposition follows. □

As a special case, we have the following:

COROLLARY 6.2. *If*

$$\left|\langle c_1(W^+), [\Sigma]\rangle\right| < n,$$

*then $\overline{\partial}_A + \overline{\partial}_A^*$ has no kernel or cokernel.*

And as another special case:

COROLLARY 6.3. *If $e = -1$, i.e.*

$$\langle c_1(W^+), [\Sigma]\rangle = -n - 2g,$$

*then the kernel of $\overline{\partial}_A + \overline{\partial}_A^*$ vanishes, and its cokernel is the $2g$-dimensional vector bundle over the Jacobian $\mathcal{J}(\Sigma) = \mathcal{M}^{\mathrm{red}}(N, W^+)$, viewed as all complex line bundles with degree $e = -1$, whose fiber over $A$ is $H^1(\Sigma, \overline{\partial}_A)$. In particular, the Chern classes of this bundle are given by*

$$c(V) = \prod_{i=1}^{g}\Big(1 + \mu(A_i)\mu(B_i)\Big),$$

*where $\{A_i, B_i\}_{i=1}^{g}$ is any symplectic basis of simple closed curves in $H_1(\Sigma; \mathbb{Z})$.*

*Proof.* Proposition 6.1 gives the identification of the cokernel bundle with the bundle over the Jacobian whose fiber over the connection $A$ (viewed as a connection in $E$ over $\Sigma$ with $\deg(E) = -1$) is $H^1(\Sigma, \overline{\partial}_A)$, a bundle which we will denote simply by $H^1(E)$. Given the identification, the statement about the Chern classes is a classical result (see [15]), but we sketch here a proof for completeness (using the Atiyah-Singer index theorem for families, see [13]).

Let $\mathbb{U}$ denote the universal complex line bundle over $\Sigma \times \mathcal{J}(\Sigma)$ parameterizing all flat bundles over $\Sigma$. Explicitly, the bundle comes equipped with a

connection
$$\nabla = d + \sum_{i=1}^{g} \Big( (\mathrm{Hol}_{A_i})A_i^* + (\mathrm{Hol}_{B_i})B_i^* \Big),$$
where the $\gamma^* \in H^1(\Sigma; \mathbb{Z})$ is the Kronecker dual to $\gamma \in H_1(\Sigma; \mathbb{Z})$, and
$$\mathrm{Hol}_\gamma \colon \mathcal{J}(\Sigma) \longrightarrow S^1$$
is the holonomy around $\gamma$. This connection parameterizes all $\overline{\partial}$ operators on the line bundles of degree zero; similarly, letting $F$ be the line bundle over $\Sigma$ of degree $2g - 1$, $\pi_\Sigma^*(F) \otimes \mathbb{U}$ parameterizes all $\overline{\partial}$ operators on the line bundles of degree $2g - 1$. Letting $d\theta$ be the volume form on $S^1$, we see that the total Chern class of $\pi_\Sigma^*(F) \otimes \mathbb{U}$ is given by
$$\begin{aligned}
c\Big(\pi_\Sigma^*(F \otimes K^{-\frac{1}{2}}) \otimes \mathbb{U}\Big) &= 1 + \pi_\Sigma^* c_1(F \otimes K^{-\frac{1}{2}}) \\
&\quad - \sum_{i=1}^{g} \Big( \mathrm{Hol}_{A_i}^*(d\theta) \cup A_i^* + \mathrm{Hol}_{B_i}^*(d\theta) \cup B_i^* \Big),
\end{aligned}$$
so its Chern character is given by
$$\begin{aligned}
\mathrm{ch}&\Big(\pi_\Sigma^*(F \otimes K^{-\frac{1}{2}}) \otimes \mathbb{U}\Big) \\
&= 1 + \pi_\Sigma^* c_1(F \otimes K^{-\frac{1}{2}}) - \sum_{i=1}^{g} \Big( \mathrm{Hol}_{A_i}^*(d\theta) \cup A_i^* + \mathrm{Hol}_{B_i}^*(d\theta) \cup B_i^* \Big) \\
&\quad + \sum_{i=1}^{g} \Big( \mathrm{Hol}_{A_i}^*(d\theta) \cup A_i^* \cup \mathrm{Hol}_{B_i}^*(d\theta) \cup B_i^* \Big).
\end{aligned}$$
Since $H^1(F) \equiv 0$, the index theorem for families gives the Chern character of the bundle $H^0(F)$ over $\mathcal{J}$ as:
$$\mathrm{ch}\Big(H^0(F)\Big) = \mathrm{ch}\Big(H^0(F) - H^1(F)\Big) = \widehat{A}(\Sigma) \mathrm{ch}\Big(\pi_\Sigma^*(F \otimes K_\Sigma^{-\frac{1}{2}}) \otimes \mathbb{U}\Big) / [\Sigma].$$
Since $\widehat{A}(\Sigma) = 1$,
$$\begin{aligned}
\mathrm{ch}&\Big(H^0(F)\Big) \\
&= \bigg(1 + c_1(F \otimes K^{-\frac{1}{2}}) - \sum \Big( \mathrm{Hol}_{A_i}^*(d\theta) \cup A_i^* + \mathrm{Hol}_{B_i}^*(d\theta) \cup B_i^* \Big) \\
&\quad + \sum_{i=1}^{g} \Big( \mathrm{Hol}_{A_i}^*(d\theta) \cup A_i^* \cup \mathrm{Hol}_{B_i}^*(d\theta) \cup B_i^* \Big) \bigg) / [\Sigma] \\
&= g - \sum_{i=1}^{g} \Big( \mathrm{Hol}_{A_i}^*(d\theta) \cup \mathrm{Hol}_{B_i}^*(d\theta) \Big).
\end{aligned}$$

Since $\mathrm{Hol}_\gamma^*(d\theta) = \mu(\gamma)$, we have that
$$c\Big(H^0(F)\Big) = \prod_{i=1}^{g} \Big(1 - \mu(A_i)\mu(B_i)\Big).$$

But $H^0(F)$ is dual to $H^1(E)$, so we have that

$$c\Big(H^1(E)\Big) = \prod_{i=1}^{g} \Big(1 + \mu(A_i)\mu(B_i)\Big),$$

as claimed. □

Note that Corollary 6.2, together with the complex interpretation of the Seiberg-Witten equations over $N$ (see the discussion surrounding Proposition 5.7) has the following consequence:

PROPOSITION 6.4. *Any solution $(A, \Phi)$ to the Seiberg-Witten equations in a $\mathrm{Spin}_{\mathbb{C}}$ structure $\mathfrak{s}$ with*

$$\Big|\langle c_1(\mathfrak{s}), [\Sigma]\rangle\Big| < n$$

*which is asymptotic to a reducible is itself also reducible.*

*Proof.* The usual integration-by-parts shows that one of $\alpha$ or $\beta$ must vanish, hence $F_A$ induces an integrable complex structure on the spinor bundle $W^+$ over $N$. Corollary 6.2 shows that the kernel of the Dirac operator for any such connection vanishes; i.e. the solution is reducible. □

## 7. Examples

Of course, symplectic surfaces of negative self-intersection in symplectic four-manifolds abound in nature, but we content ourselves here with a few (algebraic) examples.

The most obvious examples are obtained by blowing up. For instance, one could take a smooth algebraic curve of degree $d$ in $\mathbb{C}P^2$, and consider its proper transform in an $\ell$-fold blow-up $\mathbb{C}P^2 \# \ell \overline{\mathbb{C}P}^2$. This is always a genus-minimizing representative for its homology class, according Theorem 1.1.

For examples in minimal surfaces, consider the following construction. Pick three nonnegative integers $a, g, n$, and let $\mathbb{F}(n)$ denote the Hirzebruch surface of order $n$, i.e. $\mathbb{F}(n)$ is the $\mathbb{C}P^1$ bundle over $\mathbb{C}P^1$ obtained as

$$\mathbb{F}(n) = \mathbb{P}(\mathbb{C} \oplus \mathcal{O}(n)),$$

a surface with a pair of distinguished rational curves $\Sigma_+$ and $\Sigma_-$ with

$$\Sigma_\pm \cdot \Sigma_\pm = \pm n.$$

Let $C$ denote the the nodal curve obtained as the union of $2a$ distinct sections of $\mathbb{F}(n)$ (i.e. $2a$ distinct rational curves homologous to $\Sigma_+$) and $2g+2$ distinct fibers. Let $\widehat{\mathbb{F}}(n)$ denote the the surface obtained by blowing up all the crossings of $C$, and let $\widehat{C} \subset \widehat{\mathbb{F}}(n)$ denote the proper transform of $C$. Finally, let $X(a, g, n)$ denote the branched double cover of $\widehat{\mathbb{F}}(n)$ branched along $\widehat{C}$. Clearly, the preimage of $\Sigma_-$ in $X(a, g, n)$ is a curve $\Sigma$ of self-intersection number $-2n$ and genus $g$. When $g > 0$, $X(a, g, n)$ is manifestly minimal.

Nontrivial examples of Relation (4) can easily be found in blow-ups of rational or ruled surfaces. The main point here is that these manifolds have nonsimple type chambers, so, after blowing up, one gets basic classes whose pairing with the exceptional class is greater than one. Now, by representing the exceptional class by a smoothly embedded torus, we see that the hypotheses of Theorem 1.3 are satisfied, and the relation holds for nonzero invariants. In a similar vein, there are examples with any genus $g > 0$, when one uses basic classes whose pairing with the exceptional class is sufficiently large. For similar examples, see [21].

We close with some remarks concerning general properties of genus-minimizing surfaces of negative square.

The behaviour of the adjunction inequality is very sensitive to the sign of the self-intersection of a surface $\Sigma$; when the self-intersection is positive, then the genus bound for $t\Sigma$ grows with $t$ for sufficiently large $t$, while in the negative case, it drops, eventually becoming negative. This says that, in the negative case, even if the adjunction formula is sharp for $\Sigma$, it cannot remain sharp for sufficiently large multiples of $\Sigma$. On the flip side, the adjunction inequality carries extra information when it remains sharp for more than one multiple of a surface of negative square, as we shall see.

To illustrate, note that given any embedded surface $\Sigma \subset X$, for each positive integer $t$, there is an obvious *local minimizer* for $t\Sigma$, denoted $\#_t\Sigma$, which is a genus-minimizer for $t[\Sigma]$ in a neighborhood of $\Sigma$. For example, if $\Sigma$ has positive square, then by the genus-minimizing property of holomorphic curves in ruled surfaces, this genus-minimizer can be found by taking an embedded holomorphic curve in the disk bundle over $\Sigma$ representing $t[\Sigma]$; if $\Sigma$ has negative square, we take the holomorphic curve in the disk bundle, given the opposite orientation. Thus, the genus of $\#_t\Sigma$ is given by

$$g(\#_t\Sigma) = 1 + t(g-1) + \left(\frac{t(t-1)}{2}\right)|\Sigma \cdot \Sigma|.$$

Usually, for surfaces with positive square, the local minimizer in a neighborhood of a genus-minimizing surface $\Sigma$ remains globally genus-minimizing.

By contrast, if $\Sigma$ has negative square, the local minimizers for multiples of $\Sigma$ (with multiplicity greater than one) always have properly greater genus than the bound from the adjunction inequality. Thus, in cases where the adjunction formula remains sharp for more than one multiple of $\Sigma$, for example if $\Sigma$ and some $t$-fold multiple of it are both represented by symplectic submanifolds, the global genus minimizer of $t\Sigma$ has properly smaller genus than the local genus minimizer around $\Sigma$; and indeed, one can find examples where the genus minimizer of $t\Sigma$ is smaller than the genus of $\Sigma$. In the interest of concreteness, we conclude with such an example.

Fix a natural number $m$ and let $\ell = \frac{m(m-1)}{2}$. Let $L_0, \ldots, L_m \subset \mathbb{C}P^2$ be $m+1$ generic lines, generic in the sense that the intersection of any three of these lines is empty. Blowing up the $\ell$ nodes $L_i \cap L_j$ for $i \neq j \in 1, \ldots, m$, and resolving the remaining $m$ nodes $L_0 \cap L_i$ for $i = 1, \ldots, m$, we obtain a sphere in
$$X = \mathbb{C}P^2 \# \ell \overline{\mathbb{C}P^2}$$
which represents the class
$$S = (m+1)H - 2\sum_{i=1}^{\ell} E_i.$$
If $m$ is odd, this class is divisible by two; and in fact $S/2$ can be represented by a smooth, holomorphic curve of genus
$$\frac{(m-1)(m-3)}{8}.$$
(It is interesting to note that, according to Fulton [7], the nodes for $L_0 \cap L_i$ can be holomorphically resolved, to give a smooth, holomorphic representative for $S$ as well.)

In summary, we have found genus minimizers in rational surfaces with arbitrarily large genus, each of whose double is represented by an embedded sphere.

INSTITUTE FOR ADVANCED STUDY, PRINCETON, NJ
*Current address*: MICHIGAN STATE UNIVERSITY
*E-mail address*: petero@math.msu.edu

PRINCETON UNIVERSITY, PRINCETON, NJ
*E-mail address*: szabo@math.lsa.umich.edu.edu